\documentclass[10pt,a4paper]{article}
\usepackage{amsfonts}
\usepackage{amsthm}
\usepackage{amssymb}
\usepackage[centertags]{amsmath}
\usepackage{amssymb}
\usepackage{enumerate}
\usepackage{graphicx}
\providecommand{\U}[2]{\protect\rule{1.5in}{1.5in}}
\theoremstyle{plain}
\newtheorem{theorem}{Theorem}[section]

\newtheorem{definition}[theorem]{Definition}
\newtheorem{example}[theorem]{Example}
\newtheorem{lemma}[theorem]{Lemma}

\newtheorem{remarks}[theorem]{Remark}

\numberwithin{equation}{section}
\begin{document}
\title{Solving $nth$ order fuzzy differential equation by fuzzy Laplace Transform} 
\author{Latif Ahmad\footnote{a. Shaheed Benazir Bhutto University, Sheringal. b. Department of Mathematics, University of Peshawar, 25120, Khyber Pakhtunkhwa, Pakistan. E-mail: ahmad49960@yahoo.com}, Muhammad Farooq\footnote{ Department of Mathematics, University of Peshawar, 25120, Khyber Pakhtunkhwa, Pakistan. E-mail: mfarooq@upesh.edu.pk}, Saleem Abdullah\footnote{Department of Mathematics, Quaid-i-Azam University, Islamabad, Pakistan}}
\bibliographystyle{ams}
\maketitle
\begin{abstract}
In this paper, we generalize the Fuzzy Laplace Transformation (FLT) for the $nth$ derivative of a fuzzy-valued function named as $nth$ derivative theorem and under the strongly generalized differentiability concept, we use it in an analytical solution method for the solution of an $nth$ order fuzzy initial value problem (FIVP). This is a simple approach toward the solution of $nth$ order fuzzy initial value problem (FIVP) by the $nth$ generalized (FLT) form, and then we can use it to solve any order of FIVP. The related theorems and properties are proved. The method is illustrated with the help of some examples. We use MATLAB to evaluate the inverse laplace transform.
\end{abstract}
\section{Introduction}
The term fuzzy derivative has been introduced in $1972$ by Chang and Zadeh \cite{1}, while the the term fuzzy differential equation (FDE) was first formulated by Kaleva \cite{2} and Seikala \cite{3}. The theory of fuzzy differential equations (FDEs) has two branches, the one in which the Hukuhara derivative is the main tool and the other one is inclusion. This paper is based on the concept of generalized Hukuhara differentiability. The solution of FDE has a wide range of applications in the dynamic system of uncertainty. Moreover the field of FDEs is becoming a necessary part of science and real word problems \cite{2,3}. In the last few decades, many method have been applied for the solution of FIVP as discussed in \cite{4,8,20,21} but every method has advantages and disadvantages. In the near past Allahviranloo, Kaini and Barkhordari \cite{5} introduced an approach toward the existence and uniqueness to solve a second order FIVP. Here we will adopt FLT method in order to find an analytical solution of FIVP. Recently Salahshour and Allahviranloo \cite{6} has found the analytical solution of second order FIVP, then the third order FIVP has been solved by Hawrra and Amal \cite {7}. In order to solve a FIVP Allahviranloo, Kiani and Barkhordari \cite{5}  stated under what condition FLT can be applied to solve a FIVP. In \cite{5}, they proposed two conditions for the existence of solution of FIVP using FLT and its inverse, and gives some useful results in the form of first order and second derivative theorem such as linearity, continuity, uniformity and convergency under the new definition of absolute value of fuzzy-valued functions etc. They also proposed two types of absolute value of fuzzy valued function which define the convergence and exponential order of a fuzzy-valued function to find an appropriate condition. In addition they have also proved that a large class of fuzzy-valued function can be solved with the help of FLT. In this paper, we generalize FLT for $nth$ order FIVP. 
This paper is arranged as follows:\\
In section $2$, we recall some basics definitions and theorems. In section $3$, fuzzy Laplace Transform is defined. Then we prove $nth$ derivative theorem which is our contribution. In section $4$, we solve FDEs by FLT. To illustrate the method, several examples are given in section $5$. Conclusion is given in section $6$.


\section{Preliminaries}
In this section we will recall some basics definitions and theorems needed throughout the paper such as fuzzy number, fuzzy-valued function and the derivative of the fuzzy-valued functions as presented in \cite{2,3,8,9}.
\begin{definition} A fuzzy number is defined as the mapping such that $u:R\rightarrow[0,1]$, which satisfies the following four properties
\begin{enumerate}
\item $u$ is upper semi-continuous.
 \item $u$ is fuzzy  convex that is $u(\lambda  x+(1-\lambda)y) \geq \min{\{u(x), u(y)\}}. x, y\in R$ and $\lambda\in [0,1]$.
\item $u$ is normal that is $\exists$ $x_0\in R$, where $u(x_0)=1$.
\item $A=\{\overline{x \in \mathbb{R}: u(x)>0}\}$ is compact that, where $\overline{A}$ is closure of $A$.
\end{enumerate}
The mapping $f:R\rightarrow E$, where $E$ denotes the set of all fuzzy numbers(what does it mean and why it is important to write it here).
\end{definition}

\begin{definition}Obviously $R\subset{E}$ for $0\leq r \leq1$
set $[u]^{r}=\{x\in{R}:u(x)\geq r\}$ and $[u]^{0}=\{x\in{R}: u(x)>0\}$.
So it must be well known that $r\in[0,1]$ and $[u]^{r}$ is bounded closed interval.
A fuzzy number  can be defined in the parametric form as given in the following definition.(This definition does not make any sense to me)
\end{definition}
\begin{definition}
A fuzzy number in parametric form is an order pair of the form $u=(\underline{u}(r), \overline{u}(r))$, where $0\leq r\leq1$ satisfying the following conditions.
\begin{enumerate}
\item $\underline{u}(r)$ is a bounded left continuous increasing function in the interval $[0,1]$
\item $\overline{u}(r)$ is a bounded left continuous decreasing function in the interval $[0,1]$
\item $\underline{u}{(r)\leq\overline{u}(r)}$. 
\end{enumerate}
If $\underline{u}(r)=\overline{u}(r)=r$, then $r$ is called crisp number.
\end{definition}

Now we recall a triangular fuzzy number which must be in the form of
$u=(l, c, r)$ where $l,c,r\in R$ and $l\leq c\leq r$, then $\underline{u}(\alpha)=l+(c-r)\alpha$ and $\overline{u}(\alpha)=r-(r-c)\alpha$ are the end points of the $\alpha$ level set.
Since each $y\in R$ can be regarded as a fuzzy number if
\begin{eqnarray*}\widetilde{y}(t)=\begin{cases}y, \;\;\; if \;\; y=t,\\  0, \;\;\; if \;\; y\neq t.\end{cases}\end{eqnarray*}
For arbitrary fuzzy numbers $u=(\underline{u}(\alpha), \overline{u}(\alpha))$ and $v=(\underline{v}(\alpha), \overline{v}(\alpha))$ and an arbitrary crisp number $j$, we define addition and scalar multiplication as:
\begin{enumerate}
\item $(\underline{u+v})(\alpha)=(\underline{u}(\alpha)+\underline{v}(\alpha))$.
\item $(\overline{u+v})(\alpha)=(\overline{u}(\alpha)+\overline{v}(\alpha))$.
\item $(j\underline{u})(\alpha)=j\underline{u}(\alpha)$, $(j\overline{u})(\alpha)=j\overline{u}(\alpha)$ \mbox{       }  $j\geq0$
\item $(j\underline{u})(\alpha)=j\overline{u}(\alpha)\alpha, (j\overline{u})(\alpha)=j\underline{u}(\alpha)\alpha$, $j<0$
\end{enumerate}
\begin{definition} Let us suppose that x, y $\in E$, if $\exists$ $z\in E$ such that
$x=y+z$, then $z$ is called the H-difference of $x$ and $y$ and is given by $x\ominus y$\end{definition}
\begin{remarks}(see \cite{10})
Let $X$ be a cartesian product of the universes, $X_1$, $X_1, \cdots, X_n$, that is
$X=X_1 \times X_2 \times \cdots \times X_n$ and $A_{1},\cdots,A_{n}$ be $n$ fuzzy numbers in $X_1, \cdots, X_n$ respectively then  $f$ is a mapping from $X$ to a universe $Y$, and $y=f(x_{1},x_{2},\cdots,x_{n})$, then the Zadeh extension principle allows us to define a fuzzy set $B$ in $Y$ as;
 \begin{equation*}B=\{(y, u_B(y))|y=f(x_1,\cdots,x_{n}),(x_{1},\cdots,x_{n})\in X\},\end{equation*}
\noindent where
 \begin{eqnarray*}u_B(y)=\begin{cases}
 \sup_{{(x_1,\cdots,x_n)} \in f^{-1}(y)}  \min\{u_{A_1}(x_1),\cdots u_{A_n}(x_n)\}, \;\;\;  if  \;\;\; f^{-1}(y)\neq 0,\\ 0, \;\;\;\; otherwise,
 \end{cases}\end{eqnarray*}
\noindent where $f^{-1}$ is the inverse of $f$.

 The extension principle reduces in the case if $n=1$ and is given as follows: 
 $B=\{(y, u_B(y)|y=f(x), \mbox{    } x \in X)\},$
 \noindent where
 \begin{eqnarray*}u_B(y)=\begin{cases}\sup_{x\in f^{-1}(y)} \min\{u_A(x)\}, \mbox{   if   } f^{-1}(y)\neq 0,\\0, \;\;\;\; otherwise. \end{cases} \end{eqnarray*}

 By Zadeh extension principle the approximation of addition of $E$ is defined by
 $(u\oplus v)(x)=\sup_{y\in R}  \min(u(y), v(x-y))$, $x \in R$  and scalar multiplication of a fuzzy number is defined by
$(k\odot u)(x)=\{u(\frac{x}{k}), k>0$ and $0, k=0\}$, where $\widetilde{0}\in E$.
\begin{eqnarray*}(k\odot u)(x)=\begin{cases}u(\frac{x}{k}), \;\;\; k > 0,\\ 0 \;\;\; \mbox{ otherwise }.  \end{cases} \end{eqnarray*}
\end{remarks}
\noindent The Housdorff distance between the fuzzy numbers defined by \cite{6}
\[d:E\times E\longrightarrow R^{+}\cup {0},\]
\[d(u,v)=\sup_{r\in[0,1]}\max\{|\underline{u}(r)-\underline{v}(r)|, |\overline{u}(r)-\overline{v}(r)|\},\] \noindent where $u=(\underline{u}(r), \overline{u}(r))$ and $v=(\underline{v}(r), \overline{v}(r))\subset R$ has been utilized by Bede and Gal \cite{12}.(need to rewrite it and check reference \cite{12})

We know that if $d$ is a metric in $E$, then it will satisfies the following properties, introduced by Puri and Ralescu \cite{13}:
\begin{enumerate}
\item $d(u+w,v+w)=d(u,v)$, $\forall$  u, v, w $\in$ E.

\item $(k \odot u, k \odot v)=|k|d(u, v)$, $\forall$ k $\in$ R, \mbox{  and  } u, v $\in$ E.

\item $d(u \oplus v, w \oplus e)\leq d(u,w)+d(v,e)$, $\forall$ u, v, w, e $\in$  E.
\end{enumerate}
\begin{definition}(see Song and Wu \cite{14}
If $f:R\times E \longrightarrow E$, then $f$ is continuous at point $(t_0,x_0) \in R \times E$ provided that for any fixed number $r \in [0,1]$ and any $\epsilon > 0$, $\exists$ $\delta(\epsilon,r)$ such that
$d([f(t,x)]^{r}, [f(t_{0},x_{0}]^{r}) < \epsilon$
whenever $|t-t_{0}|<\delta (\epsilon, r)$ and $d([x]^{r}, [x_{0}]^{r})<\delta(\epsilon,r)$ $\forall$ t $\in$ R, x $\in E$
\end{definition}

\begin{theorem} (see Wu \cite{15})
Let $f$ be a fuzzy-valued function $[a,\infty)$ given in the parametric form as $(\underline{f}(x,\alpha), \overline{f}(x,\alpha))$ for any constant number $\alpha\in[0,1]$. Here we assume that $\underline{f}(x,\alpha)$ and $\overline{f}(x,\alpha)$ are Reman-Integral on $[a,b]$ for every $b\geq a$. Also we assume that $\underline{M}(\alpha)$ and $\overline{M}(\alpha)$ are two positive function, such that\\
$\int_a^b|\underline{f}(x,\alpha)| dx \leq \underline{M}(r)$ and $\int_a^b |\overline{f}(x,\alpha)| dx \leq \overline{M}(r)$
for every $b\geq a$, then $f(x)$ is improper integral on $[{a}, r)$. \\
Thus an improper integral will always a fuzzy number.\\
In short \[ \int_a^r f(x) dx = ( \int_a^b|\underline{f}(x,\alpha)| dx, \int_a^b |\overline{f}(x,\alpha)| dx).\]
It is will known that Hukuhare differentiability for fuzzy function was introduce by puried Ralescu in $1983$. It is based on H-differentiable.
\end{theorem}
\begin{definition}(see \cite{16})
Let $f:(a,b)\rightarrow E$ where $x_{0}\in (a,b)$, then we say that $f$ is strongly generalized differentiable at $x_0$ (Beds and Gal differentiability).
If $\exists$ an element $f'(x_0)\in E$ such that
\begin{enumerate}
 \item $\forall$ $h>0$ sufficiently small $\exists$ $f(x_0+h)\ominus f(x_0)$, $f(x_0)\ominus f(x_0-h)$, then the following limits hold (in the metric $d$)\\
 $\lim_{h\rightarrow 0}\frac{f(x_0+h)\ominus f(x_0)}{h}=\lim_{h\rightarrow 0}\frac{f(x_0)\ominus f(x_0-h)}{h}=f'(x_0)$
\noindent Or
\item $\forall h>0$ sufficiently small, $\exists$ $f(x_0)\ominus f(x_0+h)$,  $f(x_0-h)\ominus f(x_0)$, then the following limits holds (in the metric $d$)
 \item $\lim_{h\rightarrow 0}\frac{f(x_0)\ominus f(x_0+h)}{-h}=\lim_{h\rightarrow 0}\frac{f(x_0-h)\ominus f(x_0)}{-h}=f'(x_0)$
 \item $\forall$ $h>0$ sufficiently small $\exists$ $f(x_0+h)\ominus f(x_0)$,  $f(x_0-h)\ominus f(x_0)$ and the following limits holds (in metric $d$)\\
$\lim_{h\rightarrow 0}\frac{(x_0+h)\ominus f(x_0)}{h}=\lim_{h\rightarrow 0}\frac{f(x_0-h)\ominus f(x_0)}{-h}=f'(x_0)$
\item $\forall$ $h>0$ sufficiently small $\exists$ $f(x_0)\ominus f(x_0+h)$,  $f(x_0)\ominus f(x_0-h)$, then the following limits holds(in metric $d$)\\
$\lim_{h\rightarrow 0}\frac{f(x_0)\ominus f(x_0+h)}{-h}=\lim_{h\rightarrow 0}\frac{f(x_0-h)\ominus f(x_0)}{h}=f'(x_0)$
\end{enumerate}
The denominators $h$ and $-h$ denotes multiplication by $\frac{1}{h}$ $\frac{-1}{h}$ respectively.\end{definition}
\begin{theorem}(See Chalco and Reman-Flores \cite{17})
Let $F:R\rightarrow E$ be a function and denoted by $[F(t)]^{\alpha}=[f_{\alpha}(t), g_{\alpha}(t)]$ for each $\alpha\in[0,1]$, then $f$ is (i)-differentiable, then $f_{\alpha}(t)$ and $g_{\alpha}(t)$ are differentiable functions.
\end{theorem}
\begin{lemma}[see Bede and Gal\cite{18}]
Let $x_0\in R$, then the FDE $y'=f(x,y)$, $y(x_0)=y_0\in R$ and $f:R\times E\rightarrow E$ is supposed to be a continuous and is equivalent to be one of the following integral equations.
\[y(x)=y_0+\int_{x_0}^x f(t, y(t))dt \;\;\; \forall  \;\;\; x\in [x_0, x_1],\]
\noindent or
\[y(0)=y^1(x)+(-1)\odot\int_{x_0}^x f(t,y(t))dt\;\;\; \forall \;\;\; x\in [x_0, x_1],\]
\noindent on some interval $(x_0, x_1)\subset R$ depending on the strongly generalized differentiability. Integral equivalency shows that if one solution satisfies the given equation, then the other will also satisfies.
\end{lemma}
\begin{remarks}(see Gal and Bede)
In the case of strongly generalized differentiability to the FDE's $y'=f(x,y)$ we use two different integral equations. But in the case of differentiability as the definition of H-derivative, we use only one integral. The second integral equation as in Lemma $2.10$ will be in the form of $y'(t)=y'_0\ominus(-1)\int_{x_0}^x f(t,y(t))dt$.The following theorem related to the existence of solution of FIVP under the generalized differentiability (see Bede and Gal)
\end{remarks}

\begin{theorem}
Let us suppose that the following conditions are satisfied.
\begin{enumerate}
\item Let $R_0=[x_0, x_0+s]\times B(y_0, q), s,q>0, y\in E, where B(y_0,q)=\{y\in E: B(y,y_0)\leq q\}$ which denotes a closed ball in $E$ and let $f:R_0\rightarrow E$ be continuous functions such that $D(0, f(x,y))\leq M$ $\forall$ (x,y) $\in$ $R_0$ and $0\in$ E.
 \item Let $g:[x_0, x_0+s]\times [0,q]\rightarrow R$ such that $g(x, 0)\equiv 0$ and $0\leq g(x,u)\leq M$, $\forall$ x $\in [x_0, x_0+s], 0\leq u\leq q$, such that $g(x,u)$ is increasing in u, and g is such that the FIVP $u'(x)=g(x, u(x))\\ u(x)\equiv 0$ on $[x_0, x_0+s].$

\item  We have $D[f(x,y),f(x,z)\leq g(x, D(y,z))]$,$\forall$ (x,y), (x, z)$\in R_0$ and $D(y,z)\leq q.$

\item $\exists \;\; d>0$ such that for $x\in [x_0, x_0+d]$, the sequence $y'_n:[x_0, x_0+d]\rightarrow E$ given by $y'_0(x)=y_0$, $y'_n+1(x)=y_0\ominus(-1)\int_{x_0}^x f(t, y^{1}_n)dt$ defined for any $n\in N$. Then the FIVP $y'=f(x,y)$, $y(x_0)=y_0$ has two solutions that is (i)-differentiable and the other one is (ii)-differentiable for $y$.\end{enumerate}

$y^{1}=[x_0, x_0+r]\rightarrow B(y_0, q)$, where $r=\min\{s,\frac{q}{M},\frac{q}{M_1},d\}$ and the successive iteration is $y_n+1(x)=y_{0}+\int(x-{0}^{x}f(t,y^{1}_{n}(t))dt$ converges to the two solutions respectively. Now according to Theorem (2.3), we restrict our attention to function which are (i) or (ii) differentiable on their domain except on a finite number of points(see also Bede and Chalco).
\end{theorem}
\section{Fuzzy Laplace Transform}
Suppose that $f$ is a fuzzy-valued function and $p$ is a real parameter, then according to \cite{6,18} FLT of the function $f$ is defined as follows:
\begin{definition}\label{eq1}
The FLT of fuzzy-valued function is \cite{6}
\begin{equation}\label{eq2}\widehat{F}(p)=L[f(t)]=\int_{0}^{\infty}e^{-pt}f(t)dt,\end{equation}
\begin{equation}\label{eq3}\widehat{F}(p)=L[f(t)]=\lim_{\tau\rightarrow\infty}\int_{0}^{\tau}e^{-pt}f(t)dt,\end{equation}
\begin{equation}\widehat{F}(p)=[\lim_{\tau\rightarrow\infty}\int_{0}^{\tau}e^{-pt}\underline{f}(t)dt,\lim_{\tau\rightarrow\infty}\int_{0}^{\tau}e^{-pt}\overline{f}(t)dt],\end{equation}
\noindent whenever the limits exist.
\end{definition}
\begin{definition}\textbf{Classical Fuzzy Laplace Transform:} Now consider the fuzzy-valued function in which the lower and upper FLT of the function are represented by
\begin{equation}\label{eq4}\widehat{F}(p;r)=L[f(t;r)]=[l(\underline{f}(t;r)),l(\overline{f}(t;r))]\end{equation}
\noindent where
\begin{equation}\label{eq5}l[\underline{f}(t;r)]=\int_{0}^{\infty}e^{-pt}\underline{f}(t;r)dt=\lim_{\tau\rightarrow\infty} \int_{0}^{\tau}e^{-pt}\underline{f}(t;r)dt,\end{equation}
\begin{equation}\label{eq6}l[\overline{f}(t;r)]=\int_{0}^{\infty}e^{-pt}\overline{f}(t;r)dt=\lim_{\tau\rightarrow\infty}\int_{0}^{\tau}e^{-pt}\overline{f}(t;r)dt.  \end{equation}
\end{definition}
\subsection{$nth$ order fuzzy initial value problem}
In this section we are going to define an $nth$ order FIVP's under generalized H-differentiability, proposed in \cite{18}.
We define
\[y^{(n)}(t)=f(t, y(t), y'(t),y''(t),\cdots,y^{(n-1)}(t)),\]
\[y(x_0)=y_0; y'(x_0)=y'_0; y''(x_0)=y''_0,\cdots, y^{n-1}(x_0)=y^{n-1}_0,\] and
\[y(t)=(\underline{y}{(t,r)}, \overline{y}{(t,r)}),\]
\[y'(t)=(\underline{y}'(t,r), \overline{y}'(t,r)),\]
\[y''(t)=(\underline{y}''(t,r), \overline{y}''(t,r)),\]
\noindent continuing the process we get for $(n-1)th$ order that is \[y^{(n-1)}(t)=(\underline{y}^{(n-1)}(t,r), \overline{y}^{(n-1)}(t,r)),\]
\noindent are fuzzy-valued functions for $t$, where
$f(t, y(t), y'(t),y''(t),\cdots,y^{(n-1)}(t))$ is continuous fuzzy-valued function.
\begin{definition}(see \cite{8,8a}) Let $f:(a, b)\rightarrow E$ and $x_0\in (a, b)$, then the $nth$ order derivative of the function is as follows:\\
Let $f:(a,b)\rightarrow E$ where $x_{0}\in (a,b)$, then we say that $f$ is strongly generalized differentiable of the $nth$ order at $x_0$ if $\exists$ an element $f^k(x_0)\in E$ such that $\forall \;\; k=1,2\cdots,n$ one of the following holds.
\begin{enumerate}
 \item $\forall$  $h>0$ sufficiently small $\exists f^{k-1}(x_0+h)\ominus f^{k-1}(x_0)$, \\$f^{k-1}(x_0)\ominus f^{k-1}(x_0-h)$, then the following limits hold (in the metric $d$)\\
 $\lim_{h}{\rightarrow 0}\frac{f^{k-1}(x_0+h)\ominus f^{k-1}(x_0)}{h}=\lim_{h\rightarrow 0}\frac{f^{k-1}(x_0)\ominus f^{k-1}(x_0-h)}{h}=f^k(x_0)$
\noindent Or
\item $\forall$ $h>0$ sufficiently small, $\exists$ \mbox{   } $f^{k-1}(x_0)\ominus f^{k-1}(x_0+h)$,  \\$f^{k-1}(x_0-h)\ominus f^{k-1}(x_0)$, then the following limits holds (in the metric $d$)
 \item $\lim_{h\rightarrow 0}\frac{f^{k-1}(x_0)\ominus f^{k-1}(x_0+h)}{-h}=\lim_{h\rightarrow 0}\frac{f^{k-1}(x_0-h)\ominus f^{k-1}(x_0)}{-h}=f^k(x_0)$
 \item $\forall$ $h>0$ sufficiently small $\exists$ \mbox{   } $f^{k-1}(x_0+h)\ominus f^{k-1}(x_0)$, \\$f^{k-1}(x_0-h)\ominus f^{k-1}(x_0)$ and the following limits holds (in metric $d$)\\
$\lim_{h\rightarrow 0}\frac{f^{k-1}(x_0+h)\ominus f^{k-1}(x_0)}{h}=\lim_{h\rightarrow 0}\frac{f^{k-1}(x_0-h)\ominus f^{k-1}(x_0)}{-h}=f^k(x_0)$
\item $\forall$ $h>0$ sufficiently small $\exists$ \mbox{   } $f^{k-1}(x_0)\ominus f^{k-1}(x_0+h)$, \\ $f^{k-1}(x_0)\ominus f^{k-1}(x_0-h)$, then the following limits holds(in metric $d$)\\
$\lim_{h\rightarrow 0}\frac{f^{k-1}(x_0)\ominus f^{k-1}(x_0+h)}{-h}=\lim_{h\rightarrow 0}\frac{f^{k-1}(x_0-h)\ominus f^{k-1}(x_0)}{h}=f^k(x_0)$
\end{enumerate}
\end{definition}
\begin{theorem}(see \cite{7})
Let $F(t), F'(t), F''(t),\cdots,F^{(n)}(t)$ are $nth$ order differentiable fuzzy-valued functions and we denote $r$-level set of a fuzzy-valued function $F(t)$ with $[F(t)]^{r}=[f_r(t), f_r(t)]$, then $[F^{n}(t)]=[f^{n}_{r}(t), g^{n}_{r}(t)]$
\begin{proof}
Here $F(t)$ and $F'(t)$ are differentiable, then we can write as $[F''(t)]^{r}=[f''_{r}(t), g''_{r}(t)]$. Since $F''(t)$ is differentiable, then by definition $2.6$ \cite{7}, a fuzzy-valued function $F:U\rightarrow F_{0}(R^{n})$ is called Hukuhara differentiable at $t_0\in U$ if $\exists \mbox{    } DF(t_0)=F'(t_0)\in F_{0}\times R^{n}$ such that the limits
$\lim_{h\rightarrow 0}\frac{F(t_0+h)\ominus F(t_0)}{h}$ and
$\lim_{h\rightarrow 0}\frac{F(t_{0})\ominus F(t_{0}-h}{h}$ exist and is equal to $DF(t_{0})$. Similarly for $D^{2}F(t_0)$ we have
\begin{enumerate}
 \item
\begin{equation*}\begin{split}[F'(t_0+h)\ominus F'(t_0)]^{r}=[f'_{r}(t_0+h), g'_{r}(t_0+h)]\ominus[f'_{r}(t_0), g'_{r}(t_0)]\\
 =[f'_{r}(t_0+h)\ominus f'_{r}(t_0), g'_{r}(t_0+h)\ominus g'_{r}(t_0)],\end{split}\end{equation*}
\noindent and
 \item
\begin{equation*}\begin{split}[F'(t_0)\ominus F'(t_0-h)]^{r}=[f'_{r}(t_0),g'_{r}(t_0)]\ominus [f'_{r}(t_0-h), g'_{r}(t_0-h)]\\=[f'_{r}(t_0)\ominus f'_{r}(t_0-h), g'_{r}(t_0)\ominus g'_{r}(t_0-h)].\end{split}\end{equation*}
\end{enumerate}
Similarly for third order, fourth order and continuing up to $nth$ order $i.e$ $D^{n}F(t_0)$ we have
\begin{enumerate}
\item
\begin{equation*}\begin{split}[F^{(n-1)}(t_0+h)\ominus F^{(n-1)}(t_0)]^{r}=[f^{(n-1)}_{r}(t_0+h), g^{(n-1)}_{r}(t_0+h)]\ominus[f^{(n-1)}_{r}(t_0), g^{(n-1)}_{r}(t_0)]\\=[f^{(n-1)}_{r}(t_0+h)\ominus f^{(n-1)}_{r}(t_0), g^{(n-1)}_{r}(t_0+h)\ominus g^{(n-1)}_{r}(t_0)],\end{split}\end{equation*}
\begin{equation*}\begin{split}[F^{(n-1)}(t_0)\ominus F^{(n-1)}(t_0-h)]^{r}=[f^{(n-1)}_{r}(t_0),g^{(n-1)}_{r}(t_0)]\ominus [f^{(n-1)}_{r}(t_0-h), g^{(n-1)}_{r}(t_0-h)]\\=[f^{(n-1)}_{r}(t_0)\ominus f^{(n-1)}_{r}(t_0-h), g^{(n-1)}_{r}(t_0)\ominus g^{(n-1)}_{r}(t_0-h)].\end{split}\end{equation*}
\end{enumerate}
Now multiplying $\frac{1}{h}$ to the second order, third order and so on up to $nth$ order and then applying limit as $h\rightarrow 0$ on both sides we get the general form
According to \cite{19}, if $n$ is a positive integer so in the case of (1) and (2)-differentiability we can write the $nth$ derivative of the functions $F,F',\cdots ,F^{(n-1)}$ in the form of $D^{n}_{k_1\cdots k-n}F(t_0)$, where $k_i=1,2$ for $i=1,\cdots,n$ Now if we want to compute the $nth$ derivative of $F$ at $t_0$ Moreover $D^{(n-1)}_{1  1}F(t_0)$ is (1)-differentiable and $D^{(n-1)}_{2  2}F(t_0)$ is (2)-differentiable. Also $D^{(n-1)}_{1  2}F(t_0)$ is (1)and (2)-differentiable and $D^{(n-1)}_{2  1}F(t_0)$ is (2) and (1)-differentiable
 and hence proof is completed.
\end{proof}
\end{theorem}
\subsection{Convergence}
The FLT can be applied to a large number of fuzzy-valued functions \cite{6}, and in some of the examples FLT does not converge as explained below and reported in \cite{6}.
\begin{example}
Let the fuzzy-valued function $f(t)=Ce^{t^2}$, where $C\in E$, then we get
\[\lim_{\tau \rightarrow \infty} \int_{0}^{\tau}Ce^{-pt}e^{t^2}dt\rightarrow{\infty}\] for any choice of variable $p$ so the integral grows with out bounds as $\tau \rightarrow \infty $
\end{example}
In the fuzzy Laplace theory we have to use absolute value of fuzzy-valued functions. Here we will define two types of absolute value of fuzzy-valued functions as discussed in \cite{6} and is given in the following definition.
\begin{definition}
Let us consider a fuzzy-valued function whose parametric form is given in the form
\begin{equation*}f(t;r)=[\underline{f}(t;r),\overline{f}(t;r)].\end{equation*}
Now if $f$ is (1)-absolute value function, then $\forall \mbox{   } r\in [m_{1}, m_{2}]\subset [0,1]$
\begin{equation*}[f(t;r)]=[|\underline{f}(t;r)|, |\overline{f}(t;r)|].\end{equation*}
If $f$ is (2)-absolute value function, then $\forall \;r\in [m_{1}, m_{2}]\subseteq [0,1]$
\begin{equation*}[f(t;r)]=[|\overline{f}(t;r)|, |\underline{f}(t;r)|],\end{equation*}
\noindent provided that the r-cut or r-level set is satisfied by the fuzzy-valued function $|f(t;r)|$
\end{definition}
\begin{theorem}
According to \cite{6}, if a fuzzy-valued function $f$ defined as
$[f(t;r)]=[|\underline{f}(t;r)|, |\overline{f}(t;r)|]$, where $\underline{f}(t;r)$ and $\overline{f}(t;r)$ are lower and upper end points fuzzy-valued functions for $r\in [0,1]$ respectively then
\begin{enumerate}
\item If $\underline{f}(t;r)\geq 0$ $\forall r$ then $f$ is (1)-absolute value fuzzy function.
\item If $\overline{f}(t;r)\leq 0$ $\forall r$ then $f$ is (2)-absolute value fuzzy function.
\end{enumerate}
\end{theorem}
\begin{example}
Let us consider $f(t;r)=a(r)e^{t}$, \cite{6} where $a(r)=[1+r; 2-r]$, then $f(t)$ is (1)-absolute and $\forall \; r\in [0,1]$, we have
\[|f(t;r)|=[|(1+r)e^{t}|, |(2-r)e^{-t}|]=[(1+r)e^{t}, (2-r)e^{t}].\]
\end{example}
\begin{definition}
The integral (\ref{eq1}) is absolute convergent if
\begin{equation}\label{eq7}
\lim_{\tau \rightarrow \infty} \int_{0}^{\tau}|e^{-pt}f(t)|dt
\end{equation}
 \noindent exists,  that is
 \begin{equation}\label{eq8}
\lim_{\tau \rightarrow \infty}\int_{0}^{\tau}e^{-pt}|\underline{f}(t;r)|dt, \lim_{\tau \rightarrow \infty} \int_{0}^{\tau}e^{-pt}|\overline{f}(t;r)|dt
\end{equation} exist.\\\\
If $L[f(t)]$ does not converge absolutely and if $f(t)$ be (1)-absolute, then
\begin{equation}\label{eq9}
\mid\int_{\tau}^{\acute{\tau}}e^{-pt}f(t)dt\mid=[\mid\int_{\tau}^{\acute{\tau}}e^{-pt}\underline{f}(t;r)dt\mid, |\int_{\tau}^{\acute{\tau}}e^{-pt}\overline{f}(t;r)dt|],
\end{equation}
\begin{equation}\label{eq10}
|\int_{\tau}^{\acute{\tau}}e^{-pt}f(t)dt|\preceq[\int_{\tau}^{\acute{\tau}}e^{-pt}|\underline{f}(t;r)|dt,\int_{\tau}^{\acute{\tau}}e^{-pt}|\overline{f}(t;r)|dt],
\end{equation}
\begin{equation}\label{eq11}
|\int_{\tau}^{\acute{\tau}}e^{-pt}f(t)dt|=\int_{\tau}^{\acute{\tau}}e^{-pt}|f(t)|dt\rightarrow \widetilde{0},
\end{equation}
\noindent as $\tau \rightarrow {\infty}$, $\forall$ $\acute{\tau} >{\tau}$.
\noindent This implies that $L[f(t)]$ also converges. Similar case holds when $f$ is (2)-absolute.
The symbol $\leq$ is an ordering relation defined as follows:\\
\noindent For any two arbitrary fuzzy numbers $u$ and $v$, $u\leq v$ $\Leftrightarrow$ $\underline{u}(r)\leq\underline{v}(r)$ and $\overline{u}(r)\leq\overline{v}(r)$,  for all $r \in [0,1]$.
\end{definition}
Moreover in FIVP see \cite{6} in the solution of second as in definition of the integral \ref{eq1} is uniformly convergent. Also in \cite{6} the fuzzy-valued function has a jump i.e discontinuous at the point to $f$. The left and right hand limit exist but not equal. Similar also in \cite{6} the fuzzy-valued function is piece wise continuous in the interval $[0,{\infty})$ and has exponential order of $p$, then the LT of $\widehat{F}(p)=L[f(t)]$ exist for $p>s$ and converges also discussed proceeding forward. If fuzzy-valued function $f$ is piece wise continuous in $[0,{\infty})$ and has exponential order $p$, then
$\widehat{F}(p)=L[f(t)]\rightarrow 0$ as $t\rightarrow{\infty}$. For the second order FIVP \cite{6} present two theorem for translation of function as first translation and second translation theorem,
Fuzzy Laplace Transform of Derivative.
In \cite{6} Derivative Theorem.
If $f$ is continuous fuzzy-valued function on $[0,{\infty})$ then
$L[f'(t)]=pL[f(t)]\ominus f(0)$. If $f$ is (1)-differentiable. Also
$L[f'(t)]=-f(0)\ominus (-PL[f(t)])$, if $f$ is (2)-differentiable
Here,if we have to solve an $nth$ order derivative as in \cite{5} under gH-differentiability, then we will prove the results in the following equation equation$1.1$ and equation $1.2$ for $nth$ order FIVP with $n$ number of initial values.As we have proved in the previous section.
\begin{theorem}
Let suppose that  $f,f',f'',\cdots, f^{n-1}$ are continuous fuzzy-valued functions on [0,${\infty}$) and of exponential order and that $f^{(n-1)}$ is piecewise 
According to \cite{6} suppose that $f$ and $f'$ are continuous fuzzy-valued functions on $[0,\infty)$ and of exponential order and that $f''$ is piecewise continuous fuzzy-valued function on $[0,\infty)$, then
\begin{equation}\label{eq12}L(f''(t))=p^{2}L(f(t))\ominus pf(0)\ominus f'(0),
\end{equation}
\noindent if $f$ and $f'$ are (1)-differentiable
\begin{equation}\label{eq13}L(f''(t))=-f'(0)\ominus (-p^{2})L(f(t))\ominus pf(0),\end{equation}
\noindent if $f$ is (1)-differentiable and $f'$ is (2)-differentiable
\begin{equation}\label{eq14} L(f''(t))=-pf(0)\ominus (-p^{2})L(f(t))\ominus f'(0),\end{equation}
\noindent if $f$ is (2)-differentiable and $f'$ is (1)-differentiable
\begin{equation}\label{eq15}L(f''(t))=p^{2}L(f(t))\ominus pf(0)-f'(0),\end{equation}
\noindent if $f$ and $f'$ are (2)-differentiable.\\\\
\end{theorem}
\begin{theorem}
According to \cite{7}, suppose that $f(t),f'(t), f''(t)$ are the continuous fuzzy-valued function on [0,$\infty$) and of exponential order while $f'''(t)$ is piecewise continuous fuzzy-valued function on [0,$\infty$) with $f(t)=(\underline{f}(t,r),\overline{f}(t,r))$, then notations of the $nth$ derivative of the function is given by
\begin{equation}\label{eq16}L[f'''(t,r)]=p^3L[f(t)]\ominus p^2f(0)\ominus pf'(t)\ominus f''(0)\end{equation}.
\begin{proof}
Here the notation we used\\ $\underline{f}',\underline {f}'',\underline {f}'''$ \\means the lower end points functions derivatives and
\\$\overline{f}', \overline{f}'',{\overline{f}}'''$ \\are the upper end points functions derivatives, using theorem $3.2$, we have
\begin{equation}\label{eq17}
L[f'''(t)]=L[\underline{f}^{'''}(t,r), \overline{f}^{'''}(t,r)]=[l\underline{f}^{'''}(t,r), l\overline{f}^{'''}(t,r)]. \end{equation}
Now for any arbitrary fixed number $r\in [0,1]$, using the definition of classical transform, we get
\begin{equation}\label{eq18}l[\underline{f}'''(t,r)]=p^3l[\underline{f}(t,r)]-p^{2}\underline{f}(0,r)-p\underline{f}'(0,r)-\underline{f}^{''}(0,r), \end{equation}
\begin{equation}\label{eq19}l[\overline{f}^{n}(t,r)]=p^{3}l[\overline{f}(t,r)]-p^{2}\overline{f}(0,r)-p \overline{f}'(0,r)-\overline f^{''}(0,r).\end{equation}
\end{proof}
\end{theorem}
Now in order to solve the $nth$ order FIVP, we need the FLT of $nth$ derivative of the fuzzy-valued functions under the generalized H-differentiability. So we will prove the following theorem for the fuzzy Laplace Transform for $nth$ order FIVP as follows;
\begin{theorem}
Suppose that $f, f', \cdots,f^{(n-1)}$ are continuous fuzzy-valued functions on $[0, \infty)$ and of exponential order and that $f^{(n)}$ is piecewise continuous fuzzy-valued function on $[0, \infty)$, then
\begin{equation}\label{eq20}\begin{split} L(f^{(n)}(t))=p^{n}L(f(t))\ominus p^{n-1}f(0)\ominus p^{n-2}f'(0)\ominus p^{n-3}f^{''}(0)\\ \ominus\cdots \ominus f^{n-1}(0)\end{split}\end{equation}
If $f, f'\cdots f^{(n-1)}$ are (1)-differentiable
\begin{equation}\label{eq21}\begin{split} L(f^{(n)}(t))=-f^{(n-1)}(0)\ominus (-p^{n})L(f(t))\ominus p^{n-1}f(0)\ominus p^{n-2}f'(0)\\ \ominus\cdots \ominus p^{n-(n-1)}f^{(n-2)}(0) \end{split}\end{equation}
If $f, f'\cdots f^{(n-2)}$ are (1)-differentiable and $f^{(n-1)}$ is (2)-differentiable
\begin{equation}\label{eq22}\begin{split}L(f^{(n)}(t))=-p^{n-(n-1)}f^{(n-2)}\ominus f^{(n-1)}(0)\ominus (-p^{n})L(f(t))\ominus p^{n-1}f(0)\\ \ominus p^{n-2}f'(0)\ominus\cdots \ominus p^{n-(n-2)}f^{(n-3)}(0)\end{split} \end{equation}
If $f, f'\cdots f^{(n-3)}$ are (1)-differentiable and $f^{(n-1)}, f^{(n-2)}$ are (2)-differentiable
Similarly
\begin{equation}\label{eq23}\begin{split}L(f^{(n)}(t))=-p^{n-1}f(0)\ominus(-p^{n})L(f(t))\ominus p^{n-2}f'(0)\\ \ominus\cdots \ominus f^{(n-1)}(0) \end{split}\end{equation}
If $f',\cdots, f^{(n-1)}$ are (1)-differentiable and $f$ is (2)-differentiable
\\Continuing the process until we obtain $2^{n}$ system of differential equations, hence the last equation is according to \cite{18}
\begin{equation}\label{eq24}L(f^{(n)}(t))=p^{n}L(f(t))\ominus p^{n-1}f(0)\ominus p^{n-2}f'(0)\ominus p^{n-3}f''(0)\cdots -f^{n-1}(0)\end{equation}
If $f, f'\cdots f^{(n-1)}$ are (2)-differentiable
\begin{proof} \ref{eq20}
\begin{equation*}\begin{split}p^{n}L[f(t)]\ominus p^{n-1}f(0)\ominus p^{n-2}f'(0)\ominus\cdots\ominus f^{(n-1)}(0)=(p^{n}l[\underline f(t,r)]- p^{n-1}\underline f(0,r)\\-p^{n-2}\underline f'(0,r)\cdots-\underline f^{(n-1)}(0,r), p^{n}l[\overline f(t,r)]-p^{n-1}\overline f(0,r)-p^{n-2}\overline f'(0,r)\\\cdots-\overline f^{(n-1)}(0,r))\end{split}\end{equation*}
\\since
\begin{equation*}l(\underline f^{n}(t,r))=p^{n}l[\underline f(t,r)]-p^{n-1}\underline f(0,r)-p^{n-2}\underline f'(0,r)\cdots-\underline f^{(n-1)}(0,r)\end{equation*}
\begin{equation*}
l(\overline f^{(n)})=l(\overline f^{n}(t,r))=p^{n}l[\overline f(t,r)]-p^{n-1}\overline f(0,r)-p^{n-2}\overline f'(0,r)\cdots-\overline f^{(n-1)}(0,r)
\end{equation*}
where $\underline f^{(n-1)}(0,r)=\underline f^{(n-1)}(0,r)$ and $\overline f^{(n-1)}(0,r)=\overline {f^{(n-1)}}(0,r)$
\begin{equation*}\begin{split}p^{n}L[f(t)]\ominus p^{n-1}f(0)\ominus p^{n-2}f'(0)\cdots\ominus f^{(n-1)}(0)=(l(\underline f^{(n)})(t,r), l(\overline f^{(n)})(t,r))\end{split}\end{equation*}
\begin{equation*}p^{n}L[f(t)]\ominus p^{n-1}f(0)\ominus p^{n-2}f'(0)\cdots\ominus f^{(n-1)}(0)=L(\underline f^{(n)}(t,r), (\overline f^{(n)}(t,r)))\end{equation*}
\begin{equation*}p^{n}L[f(t)]\ominus p^{n-1}f(0)\ominus p^{n-2}f'(0)\cdots\ominus f^{(n-1)}(0)=L(f^{(n)}(t))\end{equation*}
Hence proof is completed
\end{proof}
Now we are going to prove the final equation (\ref{eq24}), while the middle equations are almost analogous to the proof of (\ref{eq20}) and (\ref{eq24})
\begin{proof}\ref{eq24}
\begin{equation*}\begin{split}
p^{n}L[f(t)]\ominus p^{n-1}f(0)\ominus p^{n-2}f'(0)\cdots -f^{(n-1)}(0)=(p^{n}l[\underline f(t,r)]-p^{n-1}\underline f(0,r)-\\p^{n-2}\underline f'(0,r)\cdots-\overline f^{(n-1)}(0,r), p^{n}l[\overline f(t,r)]-p^{n-1}\overline f(0,r)-p^{n-2}\overline f'(0,r)\cdots\\ -\underline f^{(n-1)}(0,r))
\end{split}\end{equation*}
since
\begin{equation*}l(\underline f^{(n)}(t,r))=l(\underline {f^{n}(t,r)})=p^{n}l[\underline f(t,r)]-p^{n-1}\underline f(0,r)- p^{n-2}\underline f'(0,r)\cdots-\underline f^{(n-1)}(0,r)\end{equation*}
\begin{equation*}
l(\overline f^{(n)})=l(\overline {f^{n}(t,r)})=p^{n}l[\overline f(t,r)]-p^{n-1}\overline f(0,r)-p^{n-2}\overline f'(0,r)\cdots-\overline f^{(n-1)}(0,r)\end{equation*}
Here we have\\
 $\overline f^{(n-1)}(0,r)=\underline f^{(n-1)}(0,r)$ and $\underline f^{(n-1)}(0,r)=\overline f^{(n-1)}$, so we know that
\begin{equation*}\begin{split}p^{n}L[f(t)]\ominus p^{n-1}f(0)\ominus p^{n-2}f'(0)\cdots -f^{(n-1)}(0)=(l(\underline f^{(n)})(t,r), l(\overline f^{(n)})(t,r))\end{split}\end{equation*}
\begin{equation*}p^{n}L[f(t)]\ominus p^{n-1}f(0)\ominus p^{n-2}f'(0)\cdots\ominus f^{(n-1)}(0)=L(\underline f^{(n)}(t,r), (\overline f^{(n)}(t,r)))\end{equation*}
\begin{equation*}p^{n}L[f(t)]\ominus p^{n-1}f(0)\ominus p^{n-2}f'(0)\cdots\ominus f^{(n-1)}(0)=L(f^{(n)}(t))\end{equation*}\\
which is the required result.
\end{proof}
\end{theorem}
\section{Constucting Solutions Via FIVP}
Consider the following $nth$ order FIVP in general form
\begin{equation}\label{eq25}y^{(n)}(t)=f(t, y(t),y'(t),\cdots, y^{(n-1)}(t)), \end{equation}
\noindent subject to the $nth$ order initial conditions $y(0)=(\underline{y}(0;r), \overline{y}(0;r))$, $y'(0)=(\underline{{y}}'(0;r), \overline{y}'(0;r))$, $y''(0)=(\underline{y}''(0;r), \overline{y}''(0;r))$.
\noindent continuing for $nth$ initial conditions
\[y^{(n-1)}(0)=(\underline{y}^{(n-1)}(0;r), \overline{y}^{(n-1)}(0;r)).\]
Now we use FLTM
\begin{equation}\label{eq26}L[y^{(n)}(t)]=L[f(t, y(t),y'(t),\cdots, y^{(n-1)}(t))]. \end{equation}
Using the theorem $5.8$ and equation (\ref{eq26})
\[p^nL[y(t)]\ominus p^{n-1}y(0)\ominus p^{n-2}y'(0)\ominus \cdots \ominus y^{(n-1)}(0)=L[f(t, y(t),y'(t),\cdots,y^{n-1}(t))].\]
The classical form
\begin{equation}\begin{split}\label{eq27}
p^{n}l[\underline{y}(t;r)]-p^{n-1}\underline{y}(0;r)-p^{n-2}\underline{y'}(0;r)-\cdots- \underline{y}^{(n-1)}(0;r)\\=l[\underline{f}(t, y(0;r),y'(0;r),\cdots,y^{(n-1)}(0;r))]
\end{split}\end{equation}
\begin{equation}\begin{split}\label{eq28}
p^{n}l[\overline{y}(t;r)]-p^{n-1}\overline{y}(0;r)-p^{n-2}\overline{y}^{'}(0;r)-\cdots- \overline{y}^{(n-1)}(0;r)\\=l[\overline{f}(t, y(0;r),y'(0;r),\cdots,y^{(n-1)}(0;r))]
\end{split}\end{equation}
In order to solve the system (\ref{eq27}) and (\ref{eq28}) we have assumed that $A(p;r)$ and $B(p;r)$ are the solutions of (\ref{eq27}) and (\ref{eq28}) respectively. So, we have
\begin{equation}\label{eq29}
l[\underline{y}(t;r)]=A(p;r),
\end{equation}
\begin{equation}\label{eq30}
l[\overline{y}(t;r)]=B(p;r).
\end{equation}
\noindent Using Inverse Laplace Transform (ILT), we have
\begin{equation}\label{eq31}
[\underline{y}(t;r)]=l^{-1}[A(p;r)],
\end{equation}
\begin{equation}\label{eq32}
[\overline{y}(t;r)]=l^{-1}[B(p;r)].
\end{equation}
\section{Examples}
\begin{example}
\begin{equation*}
y^{(iv)}=y^{'''}(t)+y^{''}(t),
\end{equation*}
\noindent with initial condition $y(0)=(3+r, 5+r)$, $y'(0)=(-3+r, -1-r)$, \\\\ $y''(0)=(8+r, 10-r)$, $y'''(0)=(12+r, 14-r)$,\\

\noindent where
\begin{equation}\label{eq33}
f(t,y(t), y'(t), y''(t), y'''(t))=y^{(iv)},\\
\end{equation}
The $r-level$ set of the FIVP is given in the following
\begin{equation}\label{eq34}
f(t,y(t),y'(t),y''(t),y'''(t))=\underline{y}'''(t;r)+\underline{y}''(t;r), \overline{y}'''(t;r)+\overline{y}''(t;r),
\end{equation}
\begin{equation}\label{eq35}
\underline{f}(t,y(t), y'(t),y''(t), y'''(t))=\underline{y}'''(t;r)+\underline{y}''(t;r),
\end{equation}
\begin{equation}\label{eq36}
\overline{f}(t,y(t), y'(t), y''(t), y'''(t))=\overline{y}'''(t;r)+\overline{y}''(t;r).
\end{equation}
Applying Laplace Transform
\begin{equation}\label{eq37}
l[\underline{f}(t, y(t), y'(t),y''(t),y'''(t))]=l[\underline{y}'''(t;r)]+l[\underline{y}''(t;r)],
\end{equation}
\begin{equation}
l[\overline{f}(t, y(t), y'(t), y''(t), y'''(t))]=l[\overline{y}'''(t;r)]+l[\overline{y}''(t;r)],
\end{equation}
\begin{equation}\begin{split}\label{eq38}
l[\underline{f}(t,y(t), y'(t),y''(t), y'''(t))]=p^3l[\underline{y}(t;r)]-p^2\underline{y}(0;r)\\-p\underline{y}'(0;r)- \underline{y}^{''}(0;r)+p^{2}l[\underline{y}(t;r)]-p\underline{y}(0;r)-\underline{y}'(0;r).
\end{split}\end{equation}
Now for upper bound, we have
\begin{equation}\begin{split}\label{eq39}
l[\overline{f}(t, y(t), y'(t), y''(t), y'''(t))]=p^3l[\overline{y}(t;r)]-p^2\overline{y}(0;r)\\-p\overline{y}'(0;r)- \overline{y}^{''}(0;r)+p^2l[\overline{y}(t;r)]-p\overline{y}(0;r)-\overline{y}'(0;r).
\end{split}\end{equation}
Now putting the initial conditions in (\ref{eq38}), we get
\begin{equation}\begin{split}\label{eq40}
l[\underline{f}(t,y(t), y'(t),y''(t), y'''(t))]=p^3l[\underline{y}(t;r)]-p^{2}(3+r)\\-p(-3+r)-(8+r)+p^{2}l[\underline{y}(t;r)]- p(3+r)-(-3+r).
\end{split}\end{equation}
Similarly putting the initial conditions in (\ref{eq39}), we have
\begin{equation}\begin{split}\label{eq41}
l[\overline{f}(t,y(t),y'(t),y''(t),y'''(t))]=p^3l[\overline{y}(t;r)]-p^2(3+r)\\-p(-3+r)-(8+r)+p^2l[\overline{y}(t;r)]-p(3+r)-(-3+r).
\end{split}\end{equation}
In general, we have
\begin{equation}\label{eq42}
L[\underline{y}^{(iv)}(t)]=L[\underline{f}(t,y(t), y'(t),y^{''}(t), y^{'''}(t))],
\end{equation}
\begin{equation}\begin{split}\label{eq43}
l[\underline{f}(t, y(t), y'(t), y^{''}(t), y^{'''}(t))]=p^4l[\underline{y}(t;r)]-p^3\underline {y}(0;r)\\-p^2\underline{y}'(0;r)-p\underline{y}^{''}(0;r)-\underline{y}^{'''}(0;r),
\end{split}\end{equation}
\begin{equation}\begin{split}\label{eq44}
[\overline{f}(t,y(t), y'(t),y^{''}(t), y^{'''}(t))]=p^4l[\overline{y}(t;r)]-p^3\overline{y}(0;r)\\-p^2\overline{y}'(0;r)- p\overline{y}^{''}(0;r)-\overline{y}^{'''}(0;r).
\end{split}\end{equation}
Now comparing equations (\ref{eq40}) and (\ref{eq43})
\begin{equation}\begin{split}\label{eq45}
p^{4}l[\underline{y}(t;r)]-p^{3}(3+r)-p^{2}(-3+r)-p(8+r)-(12+r)=p^{3}l[\underline{y}(t;r)]-\\p^2(3+r)-p(-3+r)- (8+r)+p^{2}l[\underline{y}(t;r)]-p(3+r)-(-3+r),
\end{split}\end{equation}
\begin{equation}\begin{split}\label{eq45a}
l[\underline{y}(t;r)]=(3+r)\frac{1}{p}+(-3+r)\frac{1}{p^2}\\+(8+r)\frac{p-1}{p^4-p^3-p^2}+(12+r)\frac{1}{p^4-p^3-p^2},
\end{split}\end{equation}
\begin{equation}\begin{split}\label{eq45b}
[\underline{y}(t;r)]=(3+r)l^{-1}[\frac{1}{p}]+(-3+r)l^{-1}[\frac{1}{p^2}]\\+(8+r)l^{-1}[\frac{p-1}{p^{4}-p^{3}-p^{2}}]+(12+r)l^{-1}[\frac{1}{p^{4}-p^{3}}],
\end{split}\end{equation}
\begin{equation*}\begin{split}\label{eq46}
\underline{y}(t;r)=(3+r)+(-3+r)t+(8+r)[t+2e^\frac{t}{2}\cosh (\frac{\sqrt{5}t}{2})-\frac {2}{\sqrt {5}}\sinh(\frac {\sqrt {5}t}{2})-2]\\+(12+r)[1-e^\frac {t}{2}\cos(\frac {\sqrt {5}t}{2})-\frac {3}{\sqrt {5}}\sin (\frac {\sqrt {5}t}{2})-t].
\end{split}\end{equation*}
Similarly comparing (\ref{eq41}) and (\ref{eq44}) we get
\begin{equation*}\begin{split}\label{eq47}
p^{4}l[\overline{y}(t;r)]-p^{3}(5-r)-p^{2}(-1-r)-p(10-r)-(14-r)=(p^3l[\overline{y}(t;r)]\\-p^2(5-r)-p(-1-r)- (10-r))+p^2l[\overline{y}(t;r)]-p(5-r)-(-1-r).
\end{split}\end{equation*}
After simplification we get
\begin{equation*}\label{eq48}
l[\overline{y}(t;r)]=(5-r)\frac{1}{p}+(-1-r)\frac{1}{p^2}+(10-r)\frac{p-1}{p^4-p^3-p^2}+(14-r)\frac {1}{p^4-p^3-p^2}.
\end{equation*}
\begin{equation*}\begin{split}\label{eq49}
[\overline{y}(t;r)]=(5-r)l^{-1}[\frac{1}{p}]+(-1-r)l^{-1}[\frac{1}{p^2}]+(10-r)l^{-1}[\frac{p-1}{p^4-p^3-p^2}]+(14-r)l^{-1}[\frac {1}{p^4-p^3-p^2}].
\end{split}\end{equation*}
\begin{equation*}\begin{split}\label{eq50}
\overline{y}(t;r)=(15-r)+(-1-r)t+(10-r)[t+\frac {1}{\sqrt{2}}e^\frac{t}{2}cosh(\frac{\sqrt{5}t}{2})-\frac {2}{\sqrt {5}}\sinh (\frac {\sqrt {5}t}{2})-2]+\\(14-r)[1-e^\frac {t}{2}\cos(\frac {\sqrt {5}t}{2})-\frac {3}{\sqrt {5}}\sinh (\frac {\sqrt {5}t}{2})-t].
\end{split}\end{equation*}
\end{example}
\begin{example} Solve the FIVP
\begin{equation*}\label{eq51}
y^{(iv)}=-y^{''}(t)+2y'(t)+t,\end{equation*}
\noindent subject to the initial conditions\\ $y(0)=(3-r, 1+r)$, \\$y'(0)=(4-r, -2+r)$, \\$y^{''}(0)=(7+r, 9-r)$,\\ $y^{'''}(0)=(10-r, 8+r).$
\begin{equation*}\begin{split}\label{eq52}
l[\underline{y}(t,r)]=p^4l[\underline{y}(t,r)]-p^{3}\underline{y}(0,r)-p^{2}\underline{y}'(0,r)-p\underline{y}^{''}(0,r)- \\\underline{y}^{'''}(0,r),
\end{split}\end{equation*}
\begin{equation*}\begin{split}\label{eq53}
L[\underline{y}^{''''}(t,r)]=p^4l[\underline{y}(t,r)]\ominus p^{3}\underline{y}(0,r)\ominus p^{2}\underline{y}'(0,r)\ominus \\p\underline{y}^{''}(0,r)\ominus \underline{y}^{'''}(0,r),
\end{split}\end{equation*}
\begin{equation*}\label{eq54}
l[\underline{y}(t,r)]=p^2l[\underline{y}(t,r)]-p\underline{y}(0,r)-\underline{y}'(0,r),
\end{equation*}
\begin{equation*}\label{eq55}
L[\underline{y}''(t,r)]=p^2l[\underline{y}(t,r)]-p\underline{y}(0,r)-\underline{y}'(0,r),
\end{equation*}
\begin{equation*}\label{eq56}
l[\underline{y}(t,r)]=pl[\underline{y}(t,r)]-\underline{y}(0,r),
\end{equation*}
\begin{equation*}\label{eq57}
L[\underline{y}'(t,r)]=pl[\underline{y}(t,r)]-\underline{y}(0,r).
\end{equation*}
\noindent Now putting the values, we get
\begin{equation*}\begin{split}\label{eq58}
l[\underline{y}(t,r)][p^{4}+p^{2}-2p]=p^{3}(3-r)+p^{2}(4-r)+p(7+r)\\+(10-r)+p(3-r)+(4-r)-2(3-r)+l(t),
\end{split}\end{equation*}
\begin{equation*}\begin{split}\label{eq59}
l[\underline{y}(t,r)][p^{4}+p^{2}-2p]=(p^3+p-2)(3-r)+(p^2+1)(4-r)\\+p(7+r)+(10-r+4-r-6+2r)+l(t),
\end{split}\end{equation*}
\begin{equation*}\begin{split}\label{eq60}
l[\underline{y}(t,r)][p^{4}+p^{2}-2p]=(p^3+p-2)(3-r)+(p^2+1)(4-r)+\\p(7+r)+8+l(t),
\end{split}\end{equation*}
\begin{equation*}\begin{split}\label{eq61}
\underline{y}(t,r)=(3-r)l^{-1}[\frac{p^3+p-2}{p^4+p^2-2p}]+(4-r)l^{-1}[\frac{p^2+1}{p^4+p^2-2p}]+\\(7+r)l^{-1}[\frac{p}{p^4+p^2-2p}]+8l^{-1}[\frac{1}{p^4+p^2-2p}]+6,
\end{split}\end{equation*}
\begin{equation*}\begin{split}\label{eq62}
\underline{y}(t,r)=(3-r)+(4-r)[e^\frac{t}{2}+e^\frac{t}{2}+\frac{\sqrt{7}\sin\frac{\sqrt7t}{2}}{7e^\frac{t}{2}}-1/2]+\\(7+r)[e^\frac{t}{4}-\cos(\frac{\sqrt7 t}{2})+\frac{(\frac{3\sqrt{7}\sin\frac{\sqrt{7}t}{2}}{7})}{4e^\frac{t}{2}}]+8[e^{\frac{t}{4}}\\+\cos\frac{\sqrt{7}t}{2}-\frac{(\frac{\sqrt{7}\sin\frac{\sqrt{7}t}{2}}{7})}{4e^\frac{t}{2}}-1/2]+6.
\end{split}\end{equation*}
\end{example}
\begin{example} Consider the following fourth order FIVP
\begin{equation}\label{eq63}
y^{(iv)}(t)=y^{'''}(t)-y^{''}(t),
\end{equation}
subject to the fuzzy initial conditions (FICs)\\ $y(0)=(r-1, 1-r)$,\\ $y'(0)=(r+1, 3-r)$,\\ $y''(0)=(2+r, 4-r)$,\\ $y'''(0)=(3+r,  5-r).$\\
In case of generalized H-differentiability, if we apply FLT we have a system of sixteen differential equations. It can be listed in the form of a differential operator as on the function $F(t_0)$.The list is below
\begin{table}[h]
\caption{Differentiable Operator}
\begin{center}
\begin{tabular}{|r|c|}
\hline
 1 and 2 differentiable  & 1 and 2 differentiable \\
\hline\hline
$(1)[(D_1D_1D_1D_1)F(t_0)]$ & $(9)[(D_1D_2D_1D_2)F(t_0)]$\\
 $(2)[(D_2D_1D_1D_1)F(t_0)]$ & $(10)[(D_2D_1D_2D_1)F(t_0)]$\\
 $(3)[(D_2D_2D_1D_1)F(t_0)]$ & $(11)[(D_2D_2D_1D_2)F(t_0)]$\\
$ (4)[(D_2D_2D_2D_1)F(t_0)]$& $(12)[(D_1D_2D_1D_1)F(t_0)]$\\
 $(5)[(D_2D_2D_2D_2)F(t_0)]$ & $(13)[(D_1D_1D_2D_1)F(t_0)]$\\
$(6)[(D_1D_2D_2D_2)F(t_0)]$& $(14)[(D_2D_1D_2D_2)F(t_0)]$\\
$(7)[(D_1D_1D_2D_2)F(t_0)]$& $(15)[(D_2D_1D_1D_2)F(t_0)]$\\
$(8)[(D_1D_1D_1D_2)F(t_0)]$ & $(16)[(D_1D_2D_2D_1)F(t_0)]$\\
\hline
\end{tabular}
\end{center}
\end{table}
Now here we will solve the two cases in the mentioned example that is (1)-differentiable and (2)-differentiable. First we consider that the FIVP is (1)-differentiable
Let us consider that $y(t), y'(t), \cdots ,y^{(4)}$ are (1)-differentiable
As we know that
\begin{equation*}\label{eq64}
y^{(iv)}=y'''(t)-y''(t).
\end{equation*}
Now applying FLT on both sides of the above equation, we get
\begin{equation*}\label{eq65}
L[y^{(iv)}(t)]=L[y'''(t)]-L[y''(t)],
\end{equation*}
\begin{equation*}\label{eq66}
L[f(t, y(t),y'(t),y''(t),y'''(t),y^{4}(t))]=L[y^{(iv)}],
\end{equation*}
\begin{equation*}\label{eq67}
L[y^{(iv)}]=p^4L[y(t)]\ominus p^{3}y(0)\ominus p^{2}y'(0)\ominus py''(0)\ominus y'''(0).
\end{equation*}
Now the classical FLT form of the above equation is
\begin{equation*}\begin{split}\label{eq68}
l[\underline y^{(iv)}(t,r)]=p^4l[\underline{y}(t,r)]-p^{3}\underline{y}(0,r)\ominus p^{2}\underline{y}'(0,r)- p\underline{y}^{''}(0,r)-\underline{y}^{'''}(0,r)
\end{split}\end{equation*}
\begin{equation*}\begin{split}\label{eq69}
l[\overline y^{(iv)}(t,r)]=p^4l[\overline{y}(t,r)]-p^{3}\overline{y}(0,r)-p^{2}\overline{y}'(0,r)- p\overline{y}^{''}(0,r)-\overline{y}^{'''}(0,r)
\end{split}\end{equation*}
\begin{equation*}\label{eq70}
L[y'''(t)]=p^{3}L[y(t)]\ominus p^{2}y(0)\ominus py'(0)\ominus y''(0)
\end{equation*}
The classical FLT of the above equation is
\begin{equation*}\begin{split}\label{eq71}
l[\underline y'''(t,r)]=p^{3}l[\underline{y}(t,r)]-p^{2}\underline y(0,r)-p\underline y'(0,r)-\underline y''(0,r)
\end{split}\end{equation*}
\begin{equation*}\begin{split}\label{eq72}
l[\overline y'''(t,r)]=p^{3}l[\overline {y}(t,r)]-p^{2}\overline y(0,r)-\overline y'(0,r)-\overline y''(0,r)
\end{split}\end{equation*}
\begin{equation*}\label{eq73}
L[y''(t)]=p^{2}L[y(t)]\ominus py(0)\ominus y'(0)
\end{equation*}
The classical FLT form of the above equation is
\begin{equation*}\label{eq74}
l[\underline y''(t,r)]=p^{2}l[\underline{y}(t,r)]-p\underline y(0,r)-\underline y'(0,r)
\end{equation*}
\begin{equation*}\label{eq75}
l[\overline y''(t,r)]=p^{2}l[\overline {y}(t,r)]-p\overline y(0,r)-\overline y'(0,r)
\end{equation*}
Now solve the above classical equations for lower and upper solutions, we have
\begin{equation*}\begin{split}\label{eq76}
p^4l[\underline{y}(t,r)]-p^{3}\underline{y}(0,r)-p^{2}\underline{y}'(0,r)-p\underline{y}^{''}(0,r)- \underline{y}^{'''}(0,r)=p^{3}l[\underline{y}(t,r)]-p^{2}\underline y(0,r)-\\p\underline y'(0,r)\-\underline y''(0,r)-(l[\underline y''(t,r)]=p^{2}l[\underline{y}(t,r)]\\-p\underline y(0,r)-\underline y'(0,r))
\end{split}\end{equation*}
Applying the initial conditions, we have
\begin{equation*}\begin{split}\label{eq77}
p^4l[\underline{y}(t,r)]-p^{3}(r-1)p^{2}(r+1)-p(2+r)-(3+r)=p^{3}l[\underline{y}(t,r)]-p^{2}(r-1)\\-p(r+1)-(2+r)-(p^{2}l[\underline{y}(t,r)] p(r-1)-(r+1))
\end{split}\end{equation*}
Solving the above equation for $\underline y(t,r)$, we get
\begin{equation*}\begin{split}\label{eq78}
\underline y(t,r)=(r-1)+(r+1)t+(2+r)[\frac{2}{\sqrt {3}}e^{\frac{t}{2}}\sin\frac{\sqrt {3}t}{2}-t]+(3+r)[t-e^{\frac {t}{2}}(\cos\frac {\sqrt{3}t}{2})\\+\frac {1}{\sqrt {3}}\sin\frac {\sqrt {3}t}{2}+1]
\end{split}\end{equation*}
Now we will solve the classical FLT form for $\overline y(t,r)$, we have
\begin{equation*}\begin{split}\label{eq79}
p^4l[\overline{y}(t,r)]-p^{3}(1-r)-p^{2}(3-r)-p(4-r)-(5-r)=p^{3}l[\overline {y}(t,r)]\\-p^{2}(1-r)-p(3-r)-(4-r)-(p^{2}l[\overline {y}(t,r)]-p(1-r)(3-r))
\end{split}\end{equation*}
\begin{equation*}\begin{split}\label{eq80}
\overline{y}(t,r)=(1-r)+(3-r)t+(4-r)[\frac {2}{\sqrt {3}}e^\frac {t}{2}\sin\frac {\sqrt {3}t}{2}-t]+(5-r)[t-e^\frac {t}{2}\cos\frac {\sqrt {3}t}{2}\\+\frac {1}{\sqrt {3}}\sin\frac {\sqrt {3}t}{2}+1]
\end{split}\end{equation*}
Let us consider that $y(t), y'(t), \cdots ,y^{(iv)}$ are (2)-differentiable
Now applying FLT on both sides of the above FIVP, we get
\begin{equation*}\label{eq81}
L[y^{(iv)}(t)]=L[y'''(t)]-L[y''(t)]
\end{equation*}
\begin{equation*}\label{eq82}
L[f(t, y(t),y'(t),y''(t),y'''(t),y^{iv}(t))]=L[y^{(iv)}]
\end{equation*}
Also we know that
\begin{equation*}\label{eq83}
L[y^{(iv)}]=p^4L[y(t)]\ominus p^{3}y(0)\ominus p^{2}y'(0)\ominus py''(0)\ominus y'''(0)
\end{equation*}
Now comparing the above two equations, we get the classical FLT form of equation is
\begin{equation}\begin{split}\label{eq84}
l(\underline f(t, y(t),y'(t),y''(t),y'''(t),y^{4}(t)))=p^4l[\underline{y}(t,r)]-p^{3}\underline{y}(0,r)- \\p^{2}\underline{y}'(0,r)-p\underline{y}^{''}(0,r)-\underline{y}^{'''}(0,r)\end{split} \end{equation}
\begin{equation}\begin{split}\label{eq85}
l(\overline{f}(t, y(t),y'(t),y''(t),y'''(t),y^{(iv)}(t)))=p^4l[\overline{y}(t,r)]-p^{3}\overline{y}(0,r)- \\p^{2}\overline{y}'(0,r)-\\ p\overline{y}^{''}(0,r)-\overline{y}^{'''}(0,r)\end{split} \end{equation}
In case of (2)-differentiability, we have the FIVP becomes
\begin{equation}\begin{split}\label{eq86}
l(\underline f(t, y(t),y'(t),y''(t),y'''(t),y^{(iv)}(t)))=p^{3}l[\overline {y}(t,r)]-p^{2}\overline y(0,r)\\-p\overline y'(0,r)-\overline y''(0,r)-(p^{2}l[\overline {y}(t,r)]\ominus p\overline y(0,r)-\overline y'(0,r))\end{split}\end{equation}
\begin{equation}\begin{split}\label{eq87}
l(\overline f(t, y(t),y'(t),y''(t),y'''(t),y^{(iv)}(t)))=p^{3}l[\underline{y}(t,r)]-p^{2}\underline y(0,r)\\-p\underline y'(0,r)-\underline y''(0,r)-(p^{2}l[\underline{y}(t,r)]\\-p\underline y(0,r)-\underline y'(0,r))\end{split}\end{equation}
Now comparing (\ref{eq84}), (\ref{eq85}), (\ref{eq86}) and (\ref{eq87}), we get
\begin{equation}\begin{split}\label{eq88}
p^4l[\underline{y}(t,r)]-p^{3}(r-1)-p^{2}(r+1)\-p(2+r)-(3+r)=p^{3}l[\overline {y}(t,r)]\\-p^{2}(1-r)-p(3-r)-(4-r)-(p^{2}l[\overline {y}(t,r)]-p(1-r)-(3-r))
\end{split}\end{equation}
\begin{equation}\begin{split}\label{eq89}
p^4l[\overline{y}(t,r)]-p^{3}(1-r)-p^{2}(3-r)-p(4-r)-(5-r)=p^{3}l[\underline{y}(t,r)]\\-p^{2}(r-1)-p(r+1)-(2+r)-(p^{2}l[\underline{y}(t,r)]-p(r-1)-(r+1)),
\end{split}\end{equation}
First we will solve (\ref{eq88}) and (\ref{eq89}) for $\l[\underline y(t,r)]$ and $\l[\overline y(t,r)]$
\begin{equation*}\begin{split}\label{eq90}
l[\underline{y}(t,r)]=r[\frac {p^5+p^4+2p^2+1}{p^6-p^4+2p^3-P^2}]+[\frac {3p^3+2p-4}{p^6-p^4+2p^3-P^2}],
\end{split}\end{equation*} $7$
\begin{equation*}\begin{split}\label{eq91}
l[\overline{y}(t,r)]=r[\frac {p^6-p^4+2p^3-2p^2+p-1}{p^8-p^6+2p^5-p^4}]+[\frac {3p^4-3p^3+2p^2-6p-4}{p^8-p^6+2p^5-p^4}]\end{split}.\end{equation*} $8$
Taking the inverse laplace transform, we get
\begin{equation*}\begin{split}
\underline{y}(t,r)=rl^{-1}[\frac {p^5+p^4+2p^2+1}{p^6-p^4+2p^3-P^2}]+l^{-1}[\frac {3p^3+2p-4}{p^6-p^4+2p^3-P^2}]\end{split}\end{equation*}
\begin{equation*}\begin{split}\overline{y}(t,r)=rl^{-1}[\frac {p^6-p^4+2p^3-2p^2+p-1}{p^8-p^6+2p^5-p^4}]+l^{-1}[\frac {3p^4-3p^3+2p^2-6p-4}{p^8-p^6+2p^5-p^4}]\end{split}\end{equation*}
\begin{equation*}\begin{split}\underline{y}(t,r)=r[3\cos (\frac{\sqrt {3}t}{2})+e^\frac {-t}{2}(\frac {7}{3\sqrt {5}})\sinh (\frac {\sqrt {5}t}{2})-t-2]+4t-7e^\frac {t}{2}\cos (\frac {\sqrt {3}t}{2})-\\\frac {2\sqrt {3}}{7}\sin (\frac {\sqrt {3}t}{2})-5\cosh( \frac {\sqrt {5}t}{2})-2e^\frac {t}{2}(\frac {17\sqrt {5}}{25})\sinh (\frac {\sqrt{5}t}{2})+6\end{split}\end{equation*}
\begin{equation*}\begin{split}\overline{y}(t,r)=r[3t+\frac {t^2}{2}+\frac {t^3}{6}-3\cosh(\frac {\sqrt {5}t}{2})+(\frac {7}{3\sqrt {5}})e^\frac {-t}{2}\sinh (\frac {\sqrt {5}t}{2})+3]\\+22t+7t^2+\frac {2}{3}t^3-7e^\frac {t}{2}\cos (\frac {\sqrt {3}t}{2})+(\frac {26}{7\sqrt {3}})\sin (\frac {\sqrt {3}t}{2})-59\cosh (\frac {\sqrt {5}t}{2})+(\frac {254\sqrt {5}}{295})e^\frac {t}{2}\sinh (\frac {\sqrt {5}t}{2})+33\end{split}\end{equation*}
\end{example}
\section{Conclusion}
we generalized the FLT for the $nth$ derivative of a fuzzy-valued function and provided method for the solution of an $nth$ order FIVP using the generalized differentiability concept. We have solved a number of different problems using this new approach. However some more research is needed to apply this method for the solution of system of FDEs which is in progress.
\bibliographystyle{ams}

\end{document}